\documentclass[12pt]{amsart}
\usepackage{amssymb,amscd}
\usepackage[russian]{babel}
\usepackage{graphicx}
\usepackage{pdfpages}
\usepackage{placeins}
\newtheorem{theorem}{Теорема}[section]

\newtheorem{lemma}{Лемма}[section]

\theoremstyle{definition}
\newtheorem{definition}{Определение}[section]

\theoremstyle{remark}

\ifx\pdfpagewidth\undefined \else
\fi
 \evensidemargin \oddsidemargin

\advance\hoffset by -1in
\advance\voffset by -1in

\usepackage[cp1251]{inputenc}

\author{И.~А.~Иванов-Погодаев}

\thanks{Moscow Institute of Physics and Technology}

\title{О детерминированности путей на подстановочных комплексах}

\begin{document}

\let \mathbf=\texttt
\tabcolsep 2pt

\begin{abstract}

Работа посвящена изучению комбинаторных свойств детерминированности семейства подстановочных комплексов, состоящих из четырехугольников, склееных друг с другом сторона-к-стороне. Данные свойства
являются полезными при построении алгебраических структур с конечным числом определяющих соотношений.

В частности этот метод был использован при при построении конечно определенной бесконечной нильполугруппы, удовлетворяющей тождеству $x^9=0$. Эта конструкция отвечает на проблему  Л.~Н.~Шеврина и М.~В.~Сапира. 

В данной работе исследуется возможность раскраски всего семейства комплексов в конечное число цветов, при котором выполнено свойство {\it слабой детерминированности}: если известны цвета трех вершин некоторого четырехугольника, то однозначно определен цвет четвертой стороны, кроме некоторых случаев особого расположения четырехугольника.

Даже слабой детерминированности хватает для построения конечно определенной нильполугруппы, при этом, при использовании данной конструкции, доказательство сокращается в объеме. 

Свойства детерминированности помогают корректно ввести определяющие соотношения на полугруппе путей, проходящих на построенных комплексах. Определяющие соотношения соответствуют парам эквивалентных коротких путей.

Свойства детерминированности изучались ранее в рамках теории замощений, в частности Кари и Папасоглу был построен набор   квадратных плиток, допускающий только апериодические замощения плоскости и обладающий детерминированностью: по цветам двух соседних ребер однозначно определялись цвета двух оставшихся ребер.

%Данная работа была проведена с помощью Российского Научного Фонда Грант N 17-11-01377.   Автор является победителем конкурса ``Молодая математика России''.

УДК 512.53, MSC: 20M05  

Ключевые слова: детерминированность, апериодические замощения, конечно определенные полугруппы, проблемы бернсайдовского типа; determinicity, aperiodic tilings, finitely presented semigroups, burnside-type problems.

\end{abstract}

\maketitle

\section{Введение} \label{nachalo}

В 1961 году Хао Вангом \cite{Wang} были рассмотрены квадратные плитки с разноцветными сторонами. Разные плитки можно прикладывать
друг к другу сторонами одного цвета. Был поставлен вопрос, существуют ли конечные наборы таких плиток,
с помощью которых могут быть получены только непериодические замощения плоскости. Первым такой набор был построен
Робертом Бергером \cite{Berger}, основная идея состояла в том, что замощения моделировали работу машины Тьюринга. При этом
использовалось несколько тысяч плиток. Позднее были придуманы наборы, содержащие небольшое количество плиток. Например,
интересна конструкция Рафаэля Робинсона \cite{Robinson}. Таким образом, с помошью задания конечного числа запрещенных примыканий плиток друг к другу (локальных условий) можно обеспечить глобальное свойство (апериодичность). Несколько позднее Р.Пенроузом были построены его знаменитые мозаики, а вдальнейшем в природе были открыты квазикристалы. Апериодическим замощениям посвящена обширная литература (см., например, монографию \cite{GS}.)

При построении апериодических замощений используется {\it самоподобие}. Для каждой плитки (многоугольника) некоторого конечного набора задается разбиение на меньшие плитки, являющиеся меньшими подобными копиями из того же набора. Итерация таких разбиений приводит (для нетривиальной подстановки) к апериодическим замощениям.

Одними из основных теорем теории замощений (после самого открытия апериодичности) явлились знаменитые теоремы Мозеса и Гудмана-Штраусса \cite{Mozes,C. Goodman-Strauss}. Они связывают язык локальных правил и язык подстановочных систем. Пусть задано подстановочное правило, по которому из плиток-многоугольников нескольких типов можно сложить плитку следующего уровня. Тогда стороны многоугольников могут быть декорированы конечным числом цветов так, что все замощения, соблюдающие локальные правила (плитки прикладываются друг к другу только сторонами одинаковых цветов) будут порождаться данной подстановочной системой. Таким образом, любые подстановочные системы могут быть заданы с помощью локальных правил. 

Фундаментальное значение имеет свойство {\em детерминированности} и следующим шагом явилось бы получение аналога теоремы Мозеса и Гудмана-Штраусса для разбиений с условием детерминированности. 

Для случая четырехугольных плиток детерминированность означает, что цвета двух сторон притки однозначно определяют цвета двух других сторон (существует не более одной плитки с соседними сторонами с заданной парой цветов). Это понятие возникло у Карри при исследовании групп, связанных с мозаиками (впервые такие конструкции были предложены Дж. Конвеем \cite{Conway}). Рассмотрим квадратные плитки Ванга, стороны которых раскрашены в конечное множество цветов и ставить рядом можно только плитки со стороной одинакового цвета. Существуют апериодические наборы плиток, которыми можно замостить плоскость только непериодическим образом. Как показали Кари и Папасоглу, существуют также апериодические и детерминированные наборы, где цвета двух соседних сторон однозначно определяют цвета двух других сторон. Если построить такие наборы, в которых вместе с каждой плиткой содержится и противоположная ей, то это приведет к решению проблемы Громова о существовании негиперболической группы с неположительной кривизной не содержащей $\mathbb{Z}\oplus\mathbb{Z}$.
 
Понятие {\em детерминированности} имеет значение не только в теории апериодических мозаик. Задание подстановочных замощений локальными правилами и детерминированность (в том числе, слабая) могут быть полезны при конструировании алгебраических объектов с конечным числом определяющих соотношений. Слова из полугруппы или кольца рассматриваются как пути на специально сконструированной последовательности комплексов, склеенных из 4-циклов (плиток). Пусть на вершинах и ребрах последовательности комплексов введена раскраска с глобально ограниченным числом цветов. Цветам ребер и вершин отвечают буквы конечного алфавита, вершины и ребра проходимые вдоль пути отвечают слову. Определяющим соотношениям соответствуют пары эквивалентных путей длины 2, состоящих из двух соседних сторон квадрата. Для колец и полугрупп вводятся также мономиальные соотношения для некоторых видов запрещенных путей. При этом геометрические свойства комплекса соответствуют некоторым полезным свойствам в получающемся объекте. Этот метод может быть полезен для построения конечно определенных объектов бернсайдовского типа и был применен при построении конечно определенной нильполугруппы.

Кроме того, как было указано С.~П.~Новиковым во время доклада на его семинаре, это понятие важно в теории разностных схем. При наличии детерминированности возникают операторы, позволяющие восстанавливать значение функции в четвертой вершине в трех других вершинах элемента сетки. Тем самым значение функции вдоль пути определяет значение функции в некоторой области, связанной с приклеиванием плитки. При построении ниль полугруппы раскраска вершин вдоль пути точно так же определяет участок мозаики. При этом на языке примыканий можно организовать вычислительный процесс, обеспечивающий нужный результат.
 
Таким образом, получение аналога теоремы Мозеса Гудмана-Штраусса с условием детерминированности имеет фундаметальное значение. В настоящей работе мы получим аналогичное утверждение для {\em слабой детерминированности} достаточное для построения ниль-полугруппы. 
 
 Рассмотрим четырехугольники на плоскости, разбитые на элементарные четырехугольники, примыкающие друг к другу сторона-к-стороне. Будем называть их {\em комплексами}. Пусть задано некоторое подстановочное правило, представляющее собой граф разбиения четырехугольника с пронумерованными сторонами на несколько меньших четырехугольников, стороны которых также пронумерованы. Тогда можно рассмотреть последовательность комплексов, где комплексом 1 уровня будет обычный четырехугольник, а комплекс уровня $k$ получается из комплекса уровня $k-1$ применением разбиения ко всем его элементарным четырехугольникам, с учетом ориентации сторон. Переходе от комплекса уровня $k$ к комплексу уровня $k+1$ состоит в добавлении нескольких вершин и проведении ребер. Вершины, которые добавляются на ребра комплекса уровня $k$ будем называть {\em боковыми}, а не принадлежащие ему будем называть {\em внутренними} (они добавляются внутрь каждого элементарного четырехугольника)

Рассмотрим раскраски вершин и ребер этих комплексов в конечное множество цветов. Определим понятие {\em детерминированности}, как и для квадратных плиток: раскраска последовательности комплексов называется {\em детерминированной}, если по известным цветам вершин и ребер пути из двум ребер некоторого минимального четырехугольника (плитки) однозначно восстанавливается цвет пути по другим двум ребрам этого же четырехугольника. 

\medskip

Для целей построения конечно определенных алгебраических объектов имеет смысл обощить понятие детерминированности. Можно рассматривать не только разбиения плоскости на квадратные плитки, но и более общие геометрические комплексы, склееные из топологических квадратов, когда соседние квадраты приклеиваются друг к другу сторона к стороне и одно ребро может быть стороной нескольких квадратов. Для таких комплексов также можно рассмотривать свойство детерминированности раскраски для вершин, когда цвета трех вершин одного квадрата однозначно определяют цвет четвертой вершины. Второе обобщение детерминированности для комплексов - раскрашивать не только вершины, но и ребра. Тогда можно говорить о кодировке пути на комплексе - последовательности вершин и ребер, которые проходит путь. Свойство детерминированности для путей выражается в том, что по кодировке пути по двум соседним сторонам квадрата однозначно восстанавливается кодировка пути по другим двум сторонам квадрата.

\medskip

Еще одна вариация детерминированности - слабая детерминированность,  когда такая однозначность восстановления применима не ко всем квадратам, а к почти всем, кроме некоторых специально расположенным.

Раскраска называется {\it слабо детерминированной}, если это применимо ко всем путям, кроме случая, когда четвертая вершина плитки выходит на границу некоторого подкомплекса, а три другие - нет.

Основной результат работы состоит в том, что {\em для заданной подстановочной системы всегда возможна раскраска в ограниченное число цветов, обладающая слабой детерминированностью.} Данный результат с одной стороны, нужен для построении конечно определенной ниль полугруппы, а с другой стороны, служит первым шагом для <<сильной>> детерминированности.

\medskip

\section{Определения} \label{opr}

\subsection{Структура комплексов} 

Рассмотрим четырехугольник $ABCD$, у которого зафиксированы верхняя, правая, нижняя и левая стороны.
Пусть задано разбиение $ABCD$ на конечное количество четырехугольников удовлетворяющее следующим условиям: 

1. Любые два четырехугольника либо не имеют общих точек, либо имеют общую вершину, либо имеют общую сторону;

2. Для каждого четырехугольника разбиения заданы его верхняя, правая, нижняя, левая стороны в указанном порядке по часовой стрелке.

3. На всех сторонах $ABCD$ поровну вершин четырехугольников разбиения и они разбивают стороны на одинаковые отрезки. 

4. Каждый четырехугольник разбиения примыкает к границе $0$, $1$, $2$ или $3$ вершинами, причем если примыкает двумя, то они соседние. 

\medskip

Пусть  вершины $A_1, \dots A_k$ четырехугольников разбиения расположены внутри ABCD, а вершины $U_1, \dots U_s$,
 $R_1, \dots R_s$,  $D_1, \dots D_s$,  $L_1, \dots L_s$ расположены в соответствующем порядке на верхней, правой, нижней и левой сторонах $ABCD$ соответственно. 
 
С разбиением можно связать граф, состоящий из вершин $k+4s$: $A_1, \dots A_k$, $U_1, \dots U_s$, $R_1, \dots R_s$,  $D_1, \dots D_s$,  $L_1, \dots L_s$, проведенные ребра соответствуют границам четырехугольников разбиения.

Заметим, что плиточная подстановка разбивает четырехугольник с пронумерованными сторонами на меньшие четырехугольники, также с пронумерованными сторонами.

\medskip

\begin{definition}

Четырехугольники разбиения будем называть {\it плитками} или {\it клетками}, а граф указанного выше разбиения на четырехугольники будем называть  {\it плиточной подстановкой}. Четырехугольник $ABCD$ будем называть комплексом {\it нулевого} уровня. Будем считать, что нет ребер разбиения соединяющих две вершины на периметре четырехугольника.
%Неразбитый четырехугольник будем называть комплексом нулевого уровня. 
\end{definition}

\medskip

\begin{definition}[\bf Индуктивное определение комплексов уровня $n$]
Рассмотрим четырехугольник, являющийся комплексом уровня $n$. Каждый из четырехугольников, на которые он разбит, разобьем на меньшие четырехугольники, применяя плиточную подстановку к данному четырехугольнику с пронумерованными сторонами. Получившийся подразбитый четырехугольник будем считать комплексом уровня $n+1$. 

Таким образом, плиточная подстановка $T$ определяет последовательность комплексов $K_T$.
\end{definition}

\begin{definition}
Будем называть плиточную подстановку и порождаемую ею последовательность комплексов {\it локально конечной}, если степень всех вершин всех комплексов
последовательности ограничена некоторым натуральным числом $N$.

\end{definition}

{\bf Объекты, связанные с комплексами.}

Стороны различных подкомплексов будем называть {\it макроребрами}. Ребра комплекса как графа (соединяющие соседние вершины) так и будем называть ребрами.

Рассмотрим комплекс уровня $n$.  
В нем содержатся вершины следующих видов:

1. Угловые

2. Краевые

3. Боковые.

4. Внутренние. 

Исходя из построения, комплекс уровня $n$, содержит в себе все макроребра и вершины комплексов уровней $n+1$ и более. Каждое макроребро и каждая вершина принадлежат всем комплексам, начиная с некоторого уровня. {\it Уровнем} макроребра и вершины будем называть минимальный уровень комплекса, которому они принадлежат. В частности, все макроребра на границе - нулевого уровня. Каждое ребро принадлежит некоторому макроребру, будем считать, что у него такой же уровень, как у этого макроребра.

Кроме того, каждый комплекс уровня $n$ содержит комплексы меньших уровней как подграфы. В частности, комплекс уровня $n$ есть объединение $k$ комплексов уровня $n-1$, где $k$ -- количество четырехугольников разбиения в плиточной подстановке.

\medskip

Из локальной конечности плиточной подстановки следует, что существует некоторая константа $C(N)$, такая, что если заданная вершина имеет уровень $k$, то выходящие из нее ребра будут иметь уровень не меньше, чем $k-C(N)$.

\medskip

\begin{figure}[hbtp]
\centering
\includegraphics[width=0.8\textwidth]{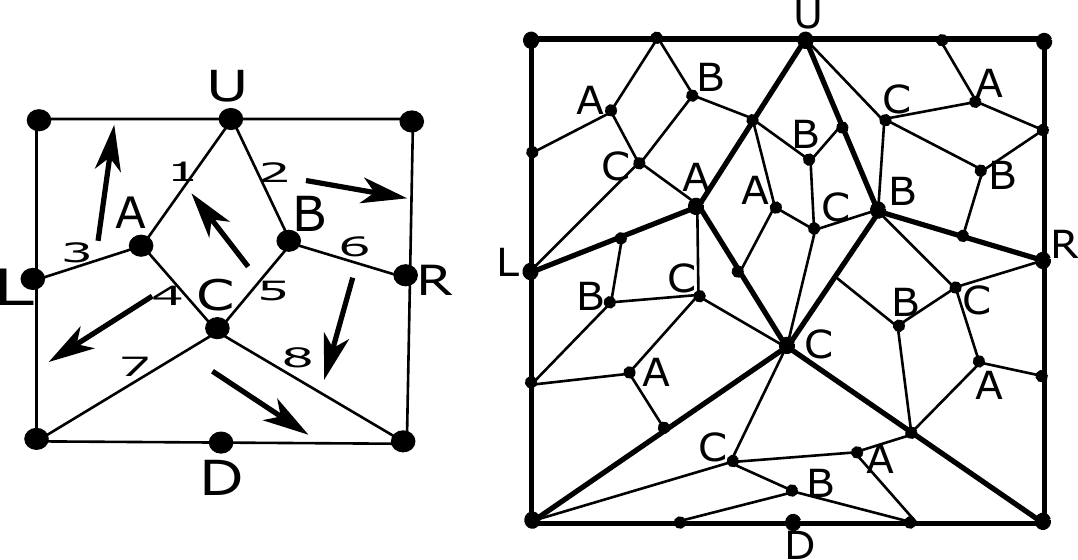}
\caption{Пример плиточной подстановки и комплекс третьего уровня с типами вершин}
\label{fig:level2}
\end{figure}

\subsection{Раскраска вершин и ребер в последовательности комплексов}

\begin{definition}[Типы макроребер] 

Рассмотрим комплекс $K_{n+1}$ уровня $n+1$. Он представляет собой объединение $k$ подкомплексов уровня $n$, где $k$ - число клеток
в плиточной подстановке. Макроребра, не лежащие на границе $K_{n+1}$, но лежащие на границах каких-то из $k$ подкомплексов, будем называть
{\it внутренними} макроребрами, принадлежащими $K_{n+1}$. 

Макроребра на границе подкомплекса, то есть такие что весь подкомплекс находится в одной из полуплоскостей 
относительно их содержащей прямой, будем называть {\it граничными} относительно него.

Пронумеруем от $1$ до $r$ все ребра в плиточной подстановке, не лежащие на границе. {\it Типами макроребер} будем называть таким образом полученные номера. Каждому внутреннему макроребру, принадлежащему комплексу будет соответствовать некоторое ребро плиточной подстановки. 
Для макроребер на границе комплекса {\it типом ребра} будем считать одну из $4$ возможных сторон (верхняя, нижняя, правая и левая). 

Для неграничных макроребер  {\it типом ребра} будем называть номер одного из внутренних ребер плиточной подстановки, соответствующего данному макроребру. Фактически тип макроребра уровня $n$ определяется, если стереть все макроребра уровней $>n$, и посмотреть, каким из внутренних ребер плиточной подстановки станет это макроребро.

Внутреннее ребро лежит на границе между двумя подкомплексами. Будем считать, что каждое внутреннеее ребро плиточной подстановки раздлеляет плоскость на две полуплоскости и они упорядочены, одна первая, другая вторая. Будем считать их сторонами этого макроребра. Каждая сторона соответствует верхней, правой, нижней или левой сторонам некоторого подкомплекса, а макроребро является общей границей этих двух подкомплексов.

\end{definition}

\medskip

\begin{definition}[Типы вершин] 

При разбиении каждой клетки комплекса $K_n$ согласно плиточной подстановке, образуются новые внутренние вершины и новые боковые и краевые вершины (внутри уже проведенных ранее макроребер. Заметим, что при дальнеших подразбиениях (переходах к комплексам больших уровней) внутренние вершины остаются внутренними, а боковые и краевые - соответственно, боковыми и краевыми. Сопоставим каждой вершине на произвольном комплексе $K_n$ {\it тип вершины}.
 
Пронумеруем внутренние вершины в плиточной подстановке от $1$ до $v$, также пронумеруем $s$ вершин на каждой стороне. Внутренним вершинам соответствуют {\it типы вершин} от $A_1$ до $A_v$. Боковая вершина лежит на некотором макроребре $l$, разделяющем два подкомплекса $U$ и $V$, одинакового уровня. Сопоставим боковой вершине $B$ принадлежащей комплексу $K_n$ ее тип. $B$ лежит на некотором макроребре, две стороны которого упорядочены. Типом боковой вершины будем считать
упорядоченную пятерку $(X,i,Y,s-i,r)$, где $X$, $Y$ -- какие-то две различные стороны из четырех (левая, верхняя, правая, нижняя), причем $X$ соответствует стороне, которой подкомплекс $U$ примыкает к макроребру $l$, а $Y$ -- стороне, которой примыкает $V$;  $i$  -- номер вершины на стороне $X$ (согласно ориентации в $U$ по часовой стрелке), $s-i$  -- номер вершины на стороне $Y$ (согласно ориентации в $V$ по часовой стрелке); $r$ - номер макроребра $l$, либо внутреннего, либо одного из четырех краевых.

\end{definition}

\begin{definition}[Уровень и скобочная ориентация] 
Для каждой вершины комплекса также определим {\it уровень}. Рассмотрим некоторый комплекс $K$. Он получается операцией разбиения из комплекса $K'$. Все добавленные при этом разбиении
вершины будем считать вершинами первого уровня. Вершины, являющиеся вершинами первого уровня для комплекса $K'$, будем называть вершинами второго уровня. Остальные вершины будем называть вершинами третьего уровня. Иными словами, только что добавленные вершины, имеют первый уровень, добавленные на предыдущем шаге -- второй, остальные -- третий.

Для боковых вершин введем скобочную ориентацию. Пусть зафиксировано некоторое макроребро и для него зафиксировано направление обхода, то есть можно говорить, какие вершины левее, какие правее. Пусть $X$ -- боковая вершина первого уровня на этом макроребре, то есть добавлена при последнем разбиении. Ровно одна из соседних с ней на данном макроребре вершина имеет второй уровень. Пусть это вершина $Y$. Будем говорить, что $X$ {\it содержит левую (открывающую) скобку}, если $Y$ стоит правее $X$. Будем говорить, что $X$ {\it содержит правую (закрывающую) скобку}, если $Y$ стоит левее $X$. 

Для вершин в середине макроребер будем говорит, что они содержат и левую и правую скобки. Таким образом, всего три возможных скобочных ориентации.

\end{definition}

\begin{definition}[Начальники вершин] 
Пусть $K_n$ подкомплекс уровня $n$. Будем называть вершины на его границе (в углах и на сторонах) его {\it начальниками}. Начальниками внутренней вершины будем называть начальником того подкомплекса, которому она принадлежит. Начальниками боковых вершин будем считать начальников того подкомплекса, которому принадлежит макроребро, на котором лежит данная вершина. Заметим, что если от $K_n$ подразбиениями его клеток перейти к комплексам больших уровней, то
начальники его вершин не поменяются. Будем считать, что все $4s+4$ начальников пронумерованы от левого верхнего угла по часовой стрелке. 
\end{definition}

\medskip

\begin{definition}[Цвет вершины] 
{\it Цветом вершины } $X$ будем называть упорядоченный набор из типа самой вершины $X$, ее уровня и скобочной ориентации.

{\it Полным цветом вершины } $X$ будем называть комбинацию цвета вершины $X$ и $4s+4$ цветов ее начальников. 
\end{definition}

\medskip
\begin{definition}[Главные ребра] 
Пусть вершина $A$ лежит на некотором макроребре. В этом случае, два выходящих из $A$ ребра тоже лежат на этом макроребре. Полуплоскости относительно макроребра упорядочены. Будем называть {\it первым главным ребром} то из них, справа от которого располагается первая из двух полуплоскостей. Второе ребро
будем называть {\it вторым главным ребром}.

\smallskip

Пусть вершина $A$ является внутренней и принадлежит подкомплексу $K$ уровня $n$. Тогда {\it главными ребрами} будем называть те выходящие из $A$ ребра,
которые тоже имеют уровень $n$ и принадлежат тому же подкомплексу $K$. Таким образом, количество главных ребер, выходящих из $A$, соответствует степени соответствующей $A$ вершины в плиточной подстановке. {\it Первым главным ребром} будем называть то ребро, которое соответствует ребру в плиточной подстановке с наименьшим номером. Заметим, что из $A$ могут выходить и другие ребра, которые будут принадлежать уже другим подкомплексам, являющимися подкомплексами $K$.

\end{definition}

\medskip

\begin{definition}[Входящие и выходящие ребра] 

Рассмотрим некоторую вершину $A$ некоторого фиксированного типа на комплексе. Смежные с $A$ ребра, и их вторые концы будем называть {\it окрестностью } вершины $A$. Заметим, что, в силу локальной конечности, существует некоторое $n$, такое что для всех вершин данного типа и уровнем выше $n$ окрестность выглядит одинаково.
В частности, если взять комплекс достаточно большого уровня, содержащий вершину $A$, то применение к нему плиточной подстановки (и переход к комплексу следующего уровня)  не добавит дополнительных ребер, смежных с $A$. Также, в силу локальной конечности, количество ребер, смежных с $A$ ограничено, причем их уровень отличается от уровня $A$ не более чем на некоторое ограниченное число. Пусть $A$ имеет уровень $n$, причем для вершин того же типа и более высокого уровня  окрестность уже не меняется.  Пронумеруем все ребра, смежные с $A$ по часовой стрелке, начиная от первого главного ребра.
\end{definition}

\medskip

{\bf  Следствия.} Из определений следуют свойства:

1.  Если выходящее и $A$ ребро $AB$ - неглавное, то как входящее в $B$ оно главное.

2.  Если выходящее и $A$ ребро $AB$ - неглавное, то вершина $A$ является одним из начальников $B$, и уровень $A$ выше, чем у $B$.

3. Если ребро $AB$ является главным входящим в $A$ и в $B$, то они лежат на одном макроребре и они имеют одинаковых начальников. При этом одна из вершин может быть внутренней.

\medskip

\begin{definition}[Кодировка пути] 
Пусть имеется путь $ABC$. Его {\it кодировкой } будем называть упорядоченный набор из полных цветов вершин $A$, $B$, $C$, и ребра входа и выхода $AB$ (относительно $A$ и $B$) и $BC$ (относительно $B$ и $C$).

\end{definition}

\medskip

\begin{figure}[hbtp]
\centering
\includegraphics[width=0.6\textwidth]{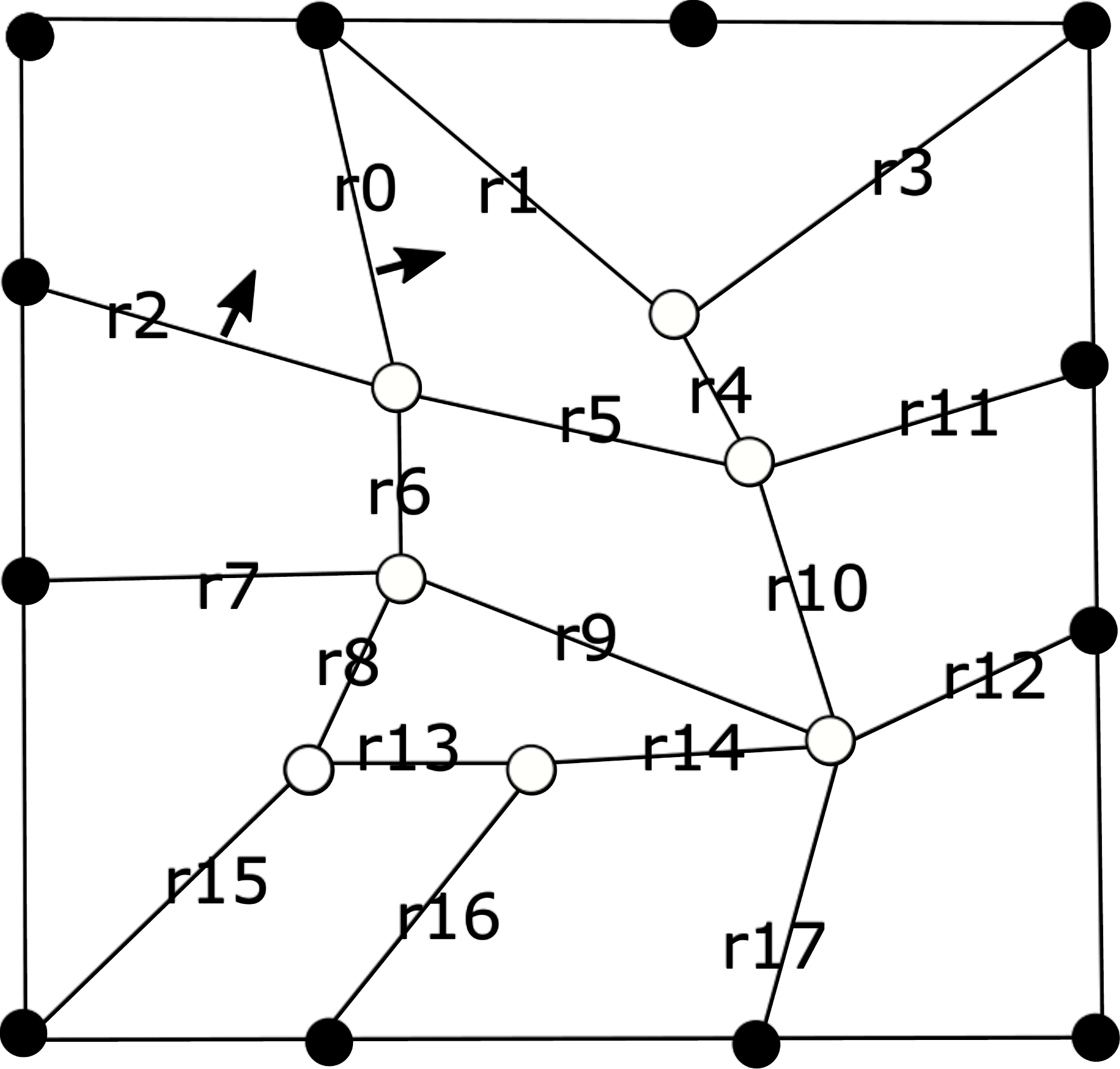}
\caption{Внутренние вершины (белые) и типы ребер внутренних ребер $r_0-r_{17}$. Стрелкой показана первая сторона (для примера, для двух ребер)  }
\label{fig:sub1}
\end{figure}

\begin{figure}[hbtp]
\centering
\includegraphics[width=0.6\textwidth]{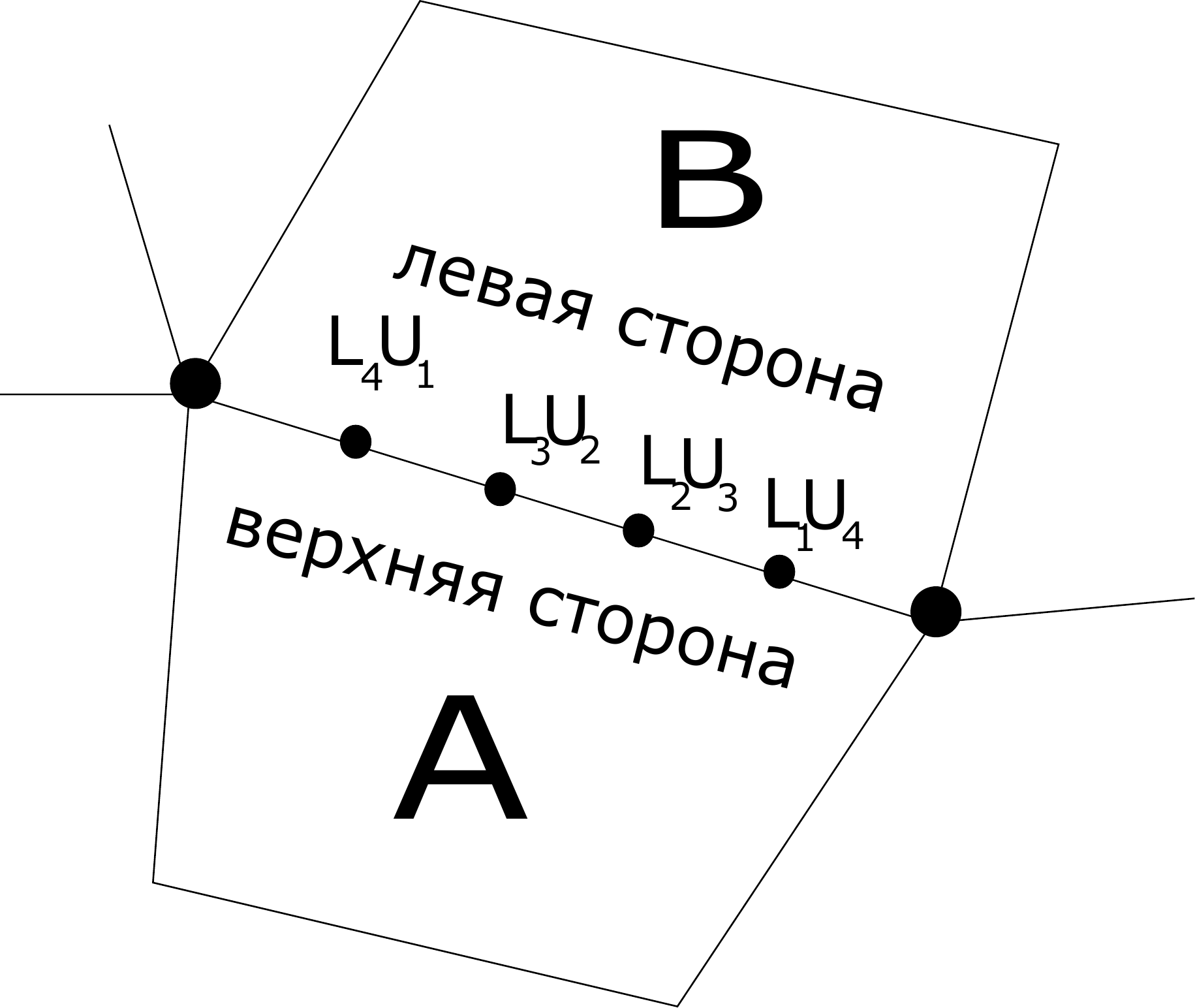}
\caption{Два подкомплекса граничат по внутреннему ребру, для одного это верхняя сторона, для другого левая, соответственно обозначаются типы вершин на внутреннем ребре, (сторона комплекса B для данного ребра первая, а другая сторона вторая) }
\label{fig:bok}
\end{figure}

\medskip

\begin{definition}
Рассмотрим некоторую вершину $X$ некоторого комплекса из последовательности, причем $X$. Пусть $X$ принадлежит элементарному четырехугольнику (клетке) $T$. Будем считать $X$ {\it особенной} относительно $T$, если ее уровень строго выше, чем уровень остальных вершин $T$. Остальные вершины будем считать неособенными.

Будем называть путь $ABC$ проходящий по двум ребрам некоторой клетки $T$ {\it регулярным}, если по его кодировке, то есть по цветам вершин
 $A$, $B$, $C$ и ребрам входа и выхода вдоль пути $ABC$ можно однозначно установить, является ли четвертая вершина $D$ клетки $T$ особенной.

Назовем раскраску {\it слабо детерминированной}, если для любого
регулярного пути $ABC$ по двум ребрам клетки $T$ и 
 неособенной четвертой вершины $D$  в той же клетке, цвет $D$, а также ребра входа и выхода вдоль пути $ADC$ однозначно определяются по 
 цветам остальных трех вершин $A$, $B$, $C$ и ребрам входа и выхода вдоль пути $ABC$ .

\end{definition}

\medskip

{\bf Комментарий.} Особенные вершины соответствуют ситуации, когда четырехугольник выходит на границу комплекса одной своей вершиной. Слабая детерминированность отличается от сильной тем, что не требуется однозначно восстанавливать цвет особенных вершин. Для приложений (в частности для построения конечно определенной нильполугруппы) достаточно слабой детерминированности по некоторым особым причинам, о которых говорится в приложении.

\medskip

\section{Слабая детерминированность} \label{raskr}

\begin{theorem}[О слабой детерминированности]
Пусть задана плиточная подстановка $\sigma$. Тогда заданная выше раскраска для $\sigma$ на последовательности комплексов $K_T$ является слабо детерминированной.

\end{theorem}

{\bf План доказательства.}  Пусть известна кодировка некоторого пути $ABC$, то есть типы вершин $A$ $B$ $C$, типы ребер $AB$, $BC$, а также типы начальников $A$, $B$ и $С$. Требуется определить, является ли данный путь двумя соседними ребрами некоторой клетки, и если является, то требуется определить тип и начальников вершины $D$ и типы ребер $AD$ и $DC$.

\medskip

\begin{definition}[О нумерации типов плиток]
Занумеруем все возможные четырехугольные клетки в плиточной подстановке. Их мы будем называть {\it типами клеток}. Таким образом, каждая плитка на комплексе любого уровня приналежит одному из этих типов.
\end{definition}

\medskip

{\bf Вспомогательные леммы.} 

\medskip

\begin{lemma}[Тип клетки определяет тип внутренних вершин] 
Пусть  известен тип клетки $T$. Тогда по нему однозначно определяются типы тех вершин в углах клетки, которые являются внутренними.  
\end{lemma}

Доказательство. Действительно, тип клетки соответствует некоторой клетке в плиточной подстановке, ее углы, если не попадают на границу плиточной подстановки,
являются внутренними вершинами, тип которых задан по плиточной подстановке.

\medskip

\begin{lemma}[Восстановление типа клетки]
Рассмотрим путь $ABC$. Пусть известен тип вершины $B$ и ребра входа в нее со стороны $A$ и $C$.  Тогда однозначно определяется тип клетки, где лежит путь
$ABC$.

\end{lemma}
Пусть известны типы вершин $A$, $B$, $C$ и ребра входов и выходов вдоль пути $ABC$. Тогда однозначно определяется, проходит ли путь $ABC$ по двум сторонам 
некоторой клетки, и ее тип.

Заметим, что для того, чтобы путь $ABC$ лежал в одной клетке необходимо и достаточно, чтобы ребра входа в $B$ были соседними.
Если тип $B$ - внутренний, то клетка однозначно определяется, как клетка, примыкающая к данной внутренней вершине, и лежащая между
 соответствующими ребрами входа.
 Аналогично, если тип $B$ - боковой или краевой, два входящих ребра полностью определяют примыкающую к ним и к $B$ клетку в плиточной подстановке.
 
 \medskip

%\begin{lemma}[Вычисление особенной ситуации]
%Рассмотрим путь $ABC$, причем его вершины лежат в одной минимальной клетке $T$. Пусть известны типы всех его %вершин и ребра входа и выхода вдоль пути. По этой информации однозначно определяется, является ли четвертая %вершина клетки особенной.
%\end{lemma}

 \begin{lemma}[Лемма об определении начальников]
 
1. Пусть $K$ является подкомплексом комплекса $T$, причем уровень $K$ отличается от $T$ на некоторую константу, зависящую от плиточной подстановки. Тогда по комбинации цветов начальников $T$ определяется комбинация цветов начальников $K$.

2. Пусть $A$ и $B$ соседние вершины на комплексе, причем $A$ является начальником $B$. Пусть известен полный цвет $A$, а также ребро выхода $AB$ и тип вершины $B$. Тогда однозначно определяется полный цвет $B$.
 
 \end{lemma}

 1. Если уровень $K$ и $T$ отличается на 1, утверждение непосредственно следует из структуры плиточной подстановки. Остальное следует по индукции по константе, так как $K$ является подкомплексом одного из комплексов, уровень которых отличается на 1 от уровня $T$.

 2. Заметим, что по ребру выхода $AB$ определяется, на сколько уровень ребра $AB$ меньше, чем уровень вершины $A$ (это следует из условия, что последовательность комплексов обладает свойством локальной конечности). При этом сам уровень $A$ мы не знаем, он может быть сколь угодно большим. Рассмотрим теперь минимальный подкомплекс $K$, содержащий ребро $AB$ в качестве внутреннего. $А$ лежит на границе $K$.
 
 Пусть $A$ принадлежит комплексу $T$. Заметим, что так как $A$ лежит на границе $K$, то $T$ не совпадает с $K$ и $K$ является его подкомплексом $T$. Тогда уровень $K$ меньше уровня $T$ на некоторую константу, зависящую от плиточной подстановки и не зависящую от уровня $T$, причем эта константа определена ребром $AB$. Тогда по первой части леммы, начальники комплекса $K$ восстанавливаются однозначно по начальникам $T$. Таким образом, начальники $B$ определяются однозначно.

\medskip

 Итак, по известному пути $ABC$ мы можем определить тип клетки, содержащую вершины $A$, $B$, $C$.

 Зная тип клетки, мы можем определить, сколько вершин из данной клетки попадает на периметр минимального содержащего ее подкомплекса (от $0$ до $3$).
  Получившиеся случаи разбираются с помощью следующих лемм. В каждом случае мы будем восстанавливать цвет оставшейся вершины $D$ и ребра входов и выходов вдоль пути $ADC$. В цвет входит тип вершины и комбинация типов начальников этой вершины.

\medskip

\begin{lemma}[Нет вершин на границе] 
Пусть $T$ -- клетка, все углы которой являются внутренними вершинами. Пусть $ABC$ -- регулярный путь, проходящий по двум сторонам $T$, и нам известна
его кодировка ( то есть цвета вершин
$A$, $B$, $C$, а также ребра входов и выходов вдоль пути $ABC$).

\begin{figure}[hbtp]
\centering
\includegraphics[width=0.3\textwidth]{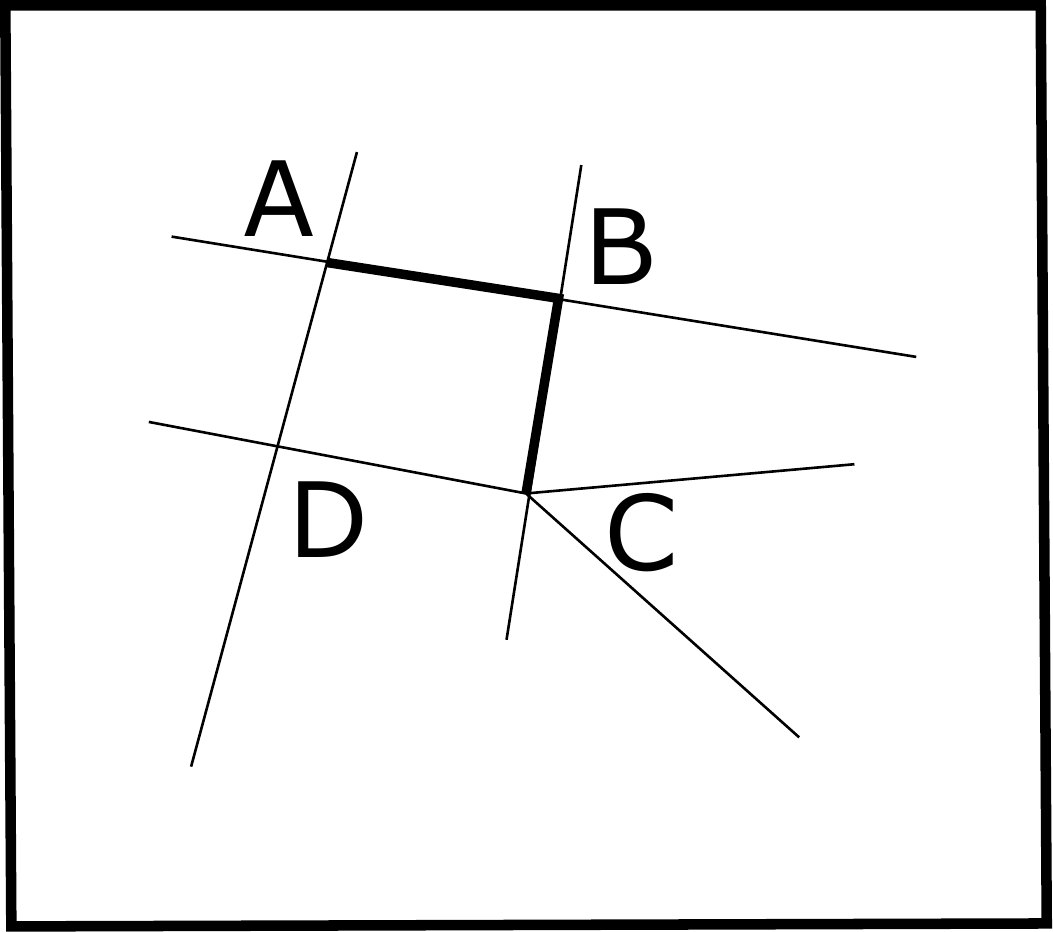}
\caption{Случай без вершин на границе}
\label{fig:case10}
\end{figure}

Тогда однозначно определяется кодировка пути $ADC$, где $D$ - четвертая вершина $T$ (то есть полный цвет $D$, ребра входов и выходов вдоль пути 
$ADC$).

\end{lemma}

Доказательство. В этом случае $T$ совпадает с одной из внутренних клеток плиточной подстановки, не имеющих общих вершин с границей, поэтому последняя вершина обязательно внутренняя и ее тип 
полностью определяется типами других вершин. Также полностью определяются ребра входов и выходов вдоль пути $ADC$,  согласно положению $T$ в плиточной подстановке. Цвета начальников $D$ такие же как у остальных вершин $A$, $B$, $C$.

\medskip

\begin{lemma}[Одна вершина на границе] 
Пусть $T$ плитка, ровно один угол которой выходит на границу плиточной подстановки (и, соответственно, на границу минимального содержащего подкомплекса). 

Пусть $ABC$ -- регулярный путь, проходящий по двум сторонам $T$, и нам известна кодировка пути $ABC$ (то есть полные цвета вершин
$A$, $B$, $C$, ребра входов и выходов вдоль пути $ABC$).

Пусть $D$ -- четвертая вершина в $T$.  Тогда либо $D$ является особенной в $T$, либо однозначно определяется кодировка пути $ADC$, где $D$ - четвертая вершина в $T$. То есть определяются ребра входов и выходов вдоль пути $ADC$, а также 
полный цвет $D$.

\begin{figure}[hbtp]
\centering
\includegraphics[width=0.3\textwidth]{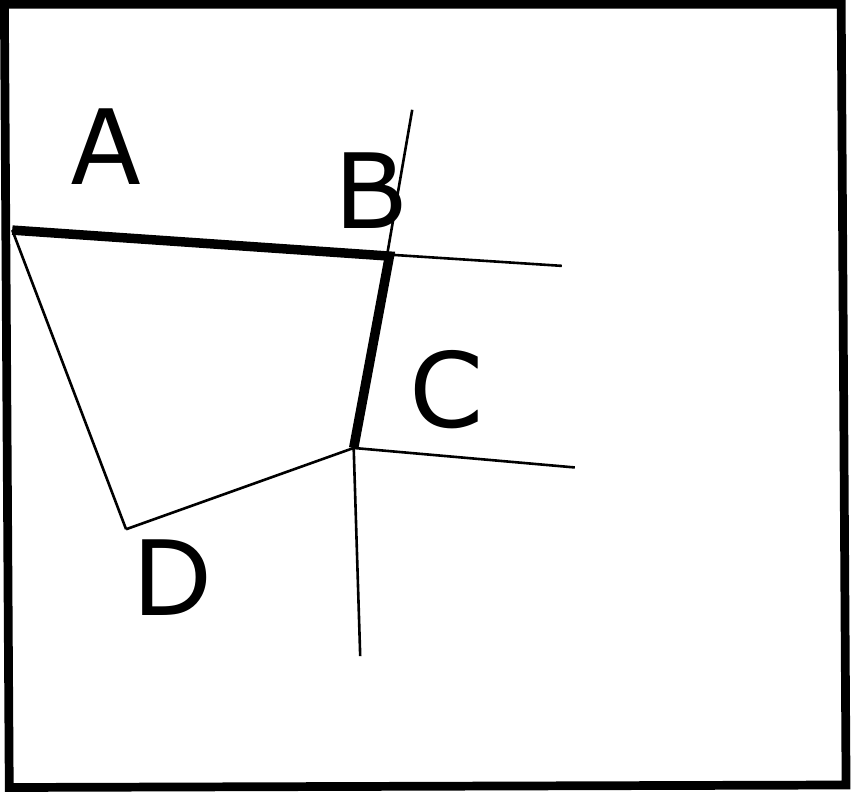}
\caption{Случай одной вершины на границе}
\label{fig:case1}
\end{figure}

\end{lemma}

Доказательство в целом аналогично случаю без вершин на границе. Если нам известны именно те три вершины, которые не лежат на границе, в этом случае 
последняя четвертая вершина является особенной. Теперь будем считать, что нам известна единственная граничная вершина $А$, и две другие вершины $B$ и $С$.

Заметим, что тип четвертой вершины $D$ однозначно восстанавливается, так как это внутренняя вершина у клетки $T$ с известным типом. Уровень $D$ первый. Кроме того, ребра входа и выхода вдоль пути $ADC$ однозначно определяются по плиточной подстановке. Комбинация начальников у $D$ такая же, как и $B$ и $C$.

%Заметим, что тогда ни одно ребро $T$ является граничным, то есть все — внутренние. Также заметим, что тип клетки полностью определяет типы трех внутренних вершин, не попадающих на границу. И зная два ребра и две вершины, мы определяем типы остальных двух ребер и оставшейся вершины.   Вторая компонента последней вершины такая же как у остальных двух.
     
  \medskip

\begin{lemma}[Две вершины на границе] 
Пусть $T$ плитка, ровно два угла которой выходят на границу плиточной подстановки (и, соответственно, на границу минимального содержащего подкомплекса). 

Пусть $ABC$ -- регулярный путь, проходящий по двум сторонам $T$, не является особенным и нам известна кодировка пути $ABC$ (то есть цвета вершин
$A$, $B$, $C$, а также ребра входов и выходов вдоль пути $ABC$).

Тогда либо четвертая вершина $D$ является особенной в плитке $T$ либо однозначно определяется кодировка пути $ADC$ (цвет $D$, ребра входов и выходов вдоль пути $ADC$).

\end{lemma}

Возможны два случая, первый, когда среди вершин $A$, $B$, $C$ есть две внутренние, и второй - когда есть одна.

\medskip

1.  Cогласно ограничениям на расположение клетки $T$ внутри плиточной подстановки, две внутренние вершины в клетке должны идти подряд. Пусть, для определенности, $B$, $C$ -- внутренние вершины. Тогда $D$ лежит на границе и является начальником $B$ и $C$ (следующим или предыдущим по часовой
стрелке после $A$, в зависимости от ориентации пути $ABC$ внутри $T$). То есть, цвет $D$ определяется. Выходящее ребро $AD$ из вершины $A$ является следующим или предыдущим выходящим ребром после $AB$. 

\medskip

Пусть выходящее из $A$ ребро $AD$ -- главное ребро. Рассмотрим это же ребро как входящее в $D$. Если это ребро $AD$ является неглавным входящим в $D$, тогда $D$, является особенной относительно $T$ (ее уровень выше, чем у $A$, а уровень $A$ выше, чем у $B$ и $С$). 

Пусть теперь $AD$ - главное ребро с обоих сторон. Тогда $D$ и $A$  имеют одинаковых начальников, то есть у $D$ они те же, что и у $A$. Итак, цвет вершины у $D$ определен, и, кроме того мы знаем входящее ребро $AD$. Кроме того, выходящее ребро $DC$ в этом случае следующее по часовой стрелке  после $DA$. То есть мы определелили все, что нужно.

\begin{figure}[hbtp]
\centering
\includegraphics[width=0.3\textwidth]{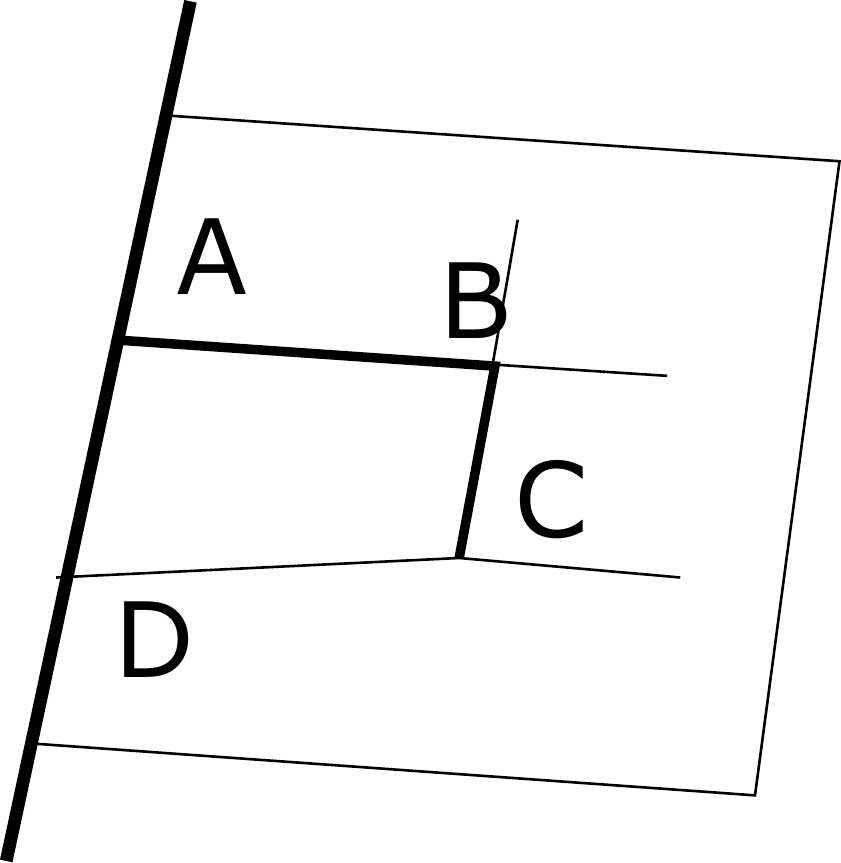}
\caption{Случай двух вершин на границе, AD - главное выходящее из А ребро}
\label{fig:case2a}
\end{figure}

Если $AD$ -- неглавное ребро, тогда $D$ лежит на неглавном ребре, выходящем из $A$. В этом случае начальники $D$ определяются по лемме определения начальников.

\begin{figure}[hbtp]
\centering
\includegraphics[width=0.3\textwidth]{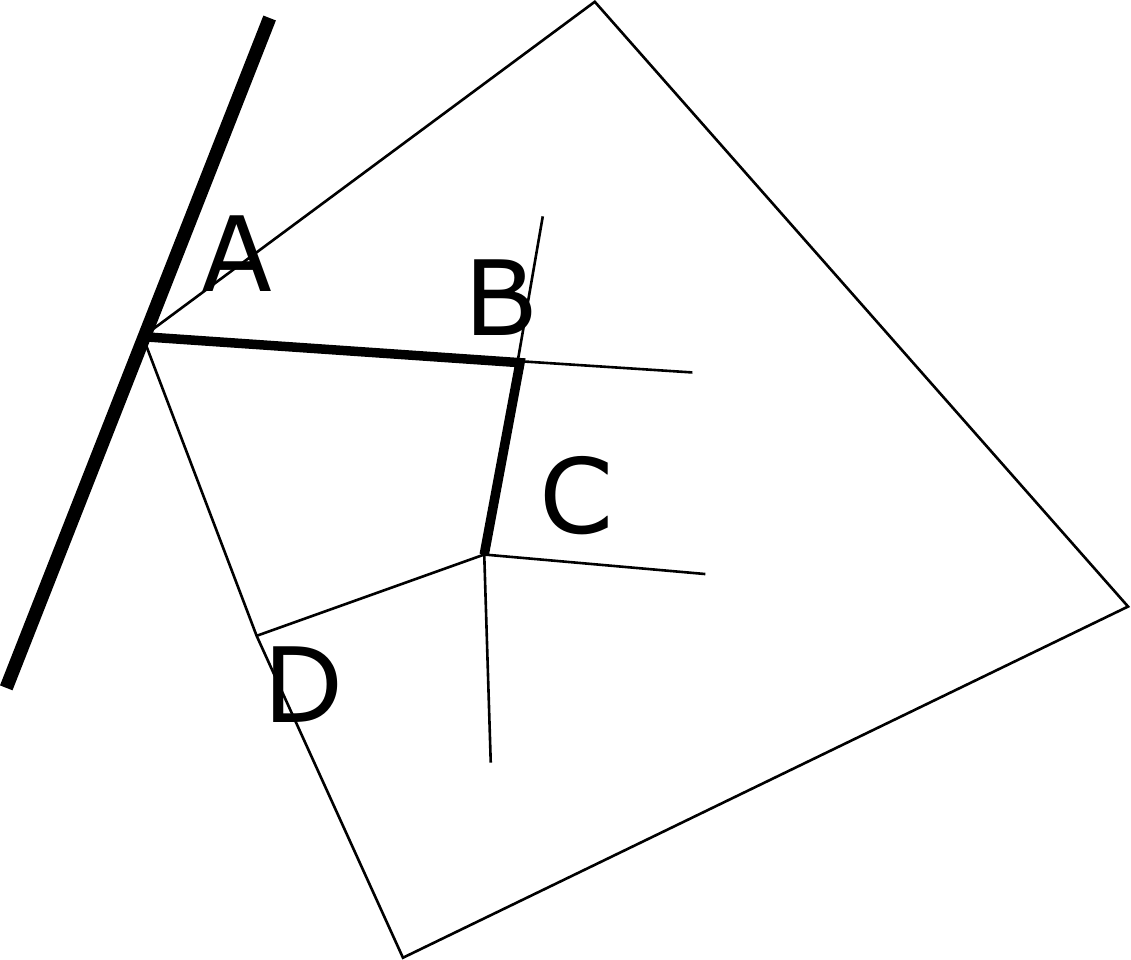}
\caption{Случай двух вершин на границе, AD - неглавное выходящее из А ребро}
\label{fig:case2b}
\end{figure}

\medskip

2. Итак, граничные вершины тоже идут подряд, пусть это $A$ и $B$, то есть $C$ и $D$ -- внутренние. Тип $D$, а также ребра входа и выхода вдоль пути $ADC$ определяются согласно расположению клетки $T$ в плиточной подстановке. Уровень $D$ первый. Комбинация цветов начальников у $D$ такая же, как у $C$.

\medskip

\begin{lemma}[Три вершины на границе]  
Пусть $T$ плитка, ровно три угла которой выходят на границу плиточной подстановки (и, соответственно, на границу минимального содержащего подкомплекса).

Пусть $ABC$ -- регулярный путь, проходящий по двум сторонам $T$, и нам известна кодировка пути $ABC$, то есть цвета вершин
$A$, $B$, $C$, а также ребра входов и выходов вдоль пути $ABC$.

Тогда, если последняя вершина $D$ не является особенной, то однозначно кодировка пути $ADC$, то есть цвет $D$ и ребра входов и выходов вдоль пути $ADC$.

\end{lemma}

Поскольку путь $ABC$ - регулярный, мы можем установить, является ли $D$ особенной. Пусть далее $D$ - не является особенной.
 
 \smallskip
 
Возможны три случая, когда вершины $A$, $B$, $C$ все граничные, когда $С$ - внутренняя  (случай с $A$ симметричен) и когда $B$ - внутренняя. 

\medskip

1. Пусть $A$, $B$, $C$ все граничные вершины. Тогда $D$ -- внутренняя и ее тип, а также ребра входа и выхода вдоль пути $ADC$ определяются согласно расположению клетки $T$ в плиточной подстановке. Уровень $D$ первый. Комбинация цветов начальников у $D$ определяется по лемме об определении начальников.

\medskip

2. Пусть $C$ - внутренняя. Выходящее ребро $AD$ из вершины $A$ тоже является следующим или предыдущим выходящим ребром после $AB$ (в зависимости от ориентации пути $ABC$ в клетке $T$). 

Пусть $AD$, как выходящее из $A$  -- главное ребро. 
 Уровень вершины $D$ не может быть выше, чем у $A$, так как в этом случае $D$ будет особенной. 
 Тогда $D$ и $A$ - вершины имеют одинаковых начальников, то есть у $D$ они те же, что и у $A$.  Кроме того, $D$ является начальником $C$, следующим или
 предыдущим перед $A$, поэтому цвет $D$ тоже определяется. Определим входящее ребро $AD$.  
 
 Пусть $D$ -боковая или краевая вершина. Входящее ребро $AD$ в этом случае -- одно из двух главных, и какое именно может быть определено по ориентации пути $ABC$.
Если $D$ - внутренняя вершина, тогда $A$ - боковая или краевая. Тогда в типе $A$ есть информация о типе макроребра, на котором лежит $A$.  По типу ребра
устанавливается откуда именно ребро $AD$ входит в $D$. Выходящее ребро $DC$ определяется как следующее или предыдущее по часовой стрелке после 
ребра $AD$ (в зависимости от ориентации $ABC$ в $T$).

\medskip
Пусть теперь $AD$ -- неглавное ребро. Тогда $D$ лежит на неглавном ребре, выходящем из $A$, $D$ - боковая или краевая вершина.  $D$ является начальником $C$, следующим или
 предыдущим перед $A$, поэтому цвет $D$ определяется. Начальники $D$ определяются по лемме об определении начальников. Заметим, что
 как входящее в $D$ ребро $AD$ -- главное и какое именно из двух, определяется согласно ориентации $ABC$ в $T$. Выходящее ребро $DC$ - следующее или 
 предыдущее после $AD$.

\medskip

3.  Пусть $B$ - внутренняя. Полностью аналогичный случаю 2 разбор в зависимости от того, главное ли ребро $AD$.

\medskip

Из приведенных выше лемм следует доказательство теоремы о слабой детерминированности.

\medskip

{\bf Следствие из слабой детерминированности.} Пусть плиточная подстановка представляет собой квадрат $n\times n$. В этом случае комплексы $K_n$ будут предствлять собой конечные подквадраты $n^k\times n^k$ квадратной решетки. Пусть задана слабодетерминированная раскраска. Выпишем все кодировки $AuBvC$ пути длины $4$ на получившейся последовательности комплексов. Теперь объявим запрещенными все кодировки, не встречающиеся на последовательности комплексов $K_n$. Рассмотрим теперь прозвольный квадрат $T$ на квадратной решетке, вершины и ребра которой раскрашены в определенное нами множество цветов, причем правильным образом. То есть на квадрате не встречается запрещенных путей. Тогда $T$ принадлежит некоторому комплексу $K_m$.

\smallskip

{\bf Макет доказательства}
Для квадратов $2 \times 2$ утверждение очевидно. Пусть утверждение доказано для квадратов, менее $2m$. Разобьем квадрат $2m \times m$ на два одинаковых квадрата. Каждый из них принадлежит некоторому комплексу по предположению индукции. Тогда рассматривая границу и используя детерминированность, можно <<сшить>> два квадрата и получать комплекс, содержащий оба квадрата сразу.

\medskip

{\bf Замечание.}  Данное следствие представляет собой вариацию теоремы Гудмана-Штраусса для квадратной решетки. Использую слабую детерминированность, похожим способом можно доказать эту теорему и для произвавольных комплексов.

\medskip

\section{Приложение} \label{appendix}

{\bf О цели введения понятия слабой детеминированности.}  При построении конечно определенной полугруппы на последовательности комплексов свойство детерминированности используется чтобы корректно ввести определяющие соотношения. То есть для каждой пары путей, с общими концами, обходящими плитку с разных сторон вводится соотношение, приравнивающее их кодировки. 

Пути, для которых слабая детерминированность отличается от обычной, соответствуют обходу плитки, одним углом выходящией на границу плиточной подстановки. Будем называть такие пути особыми. (То есть особый путь нельзя перекинуть на другую сторону плитки, так как при этом неодозначно определяется цвет вершины на границе, которая не входит в особый путь, но входит в парный к нему). При введении определяющих соотношений мы также приравниваем к нулю короткие некратчайшие пути. Для некоторых видов плиточных подстановок удается доказать, что если достаточно длинный путь включает особый подпуть, то его можно локально преобразовать к виду, содержащему некратчайший участок. То есть, такой путь будет все равно нулевой. Это дает возможность не вводить определяющих соотношений для особых путей, и обойтись слабой детерминированностью.

\medskip

В этом приложении свойство слабой детерминированности расширяется на класс комплексов с подклейками, которые используются для построения конечно определенной нильполугруппы.

\medskip

Сначала определим последовательность комплексов с подклейками. Комплекс уровня $n$ будет получаться из комплекса предыдущего уровня посредством применения {\it двух} операций: помимо разбиения, будет также проводиться подклейка - добавление в комплекс плитки, три угла которой лежат в комплексе, а четвертый угол является дополнительной новой вершиной.

Пусть задана плиточная подстановка $\sigma$. Определим последовательность комплексов, связанных с этой подстановкой.
Будем считать, что обычный квадрат является комплексом уровня $0$. Граф, совпадающий с плиточной подстановкой $\sigma$, будем считать комплексом уровня $1$. Заметим, что каждый четырехугольник разбиения ориентирован (для него определены верхняя, левая, правая и нижняя стороны) и к нему можно применить плиточную подстановку, подразбивающую его на более мелкие четырехугольники. Получившийся граф будем считать 
комплексом уровня $2$.

Комплексу уровня $0$ принадлежат углы, будем считать, что они имеют уровень $0$. Будем также считать, что вершины комплекса уровня $2$, ему {\it принадлежат}, если они не являются вершинами комплекса уровня $1$, аналогично для комплекса уровня $1$ и $0$. Если вершина принадлежит комплексу уровня $k$, то она имеет уровень $k$.

\medskip

\begin{definition}
Пусть $K_n$ -- комплекс уровня $n$. Рассмотрим пути $ABC$, такие что $B$ имеет уровень $n-1$, а вершины $A$ и $C$ имеют уровень $n$. Пусть также, в комплексе не существует клетки, содержащей вершины $A$, $B$, $C$ в качестве своих углов.  Для каждого такого пути введем новую вершину $D$ и добавим ребра $AD$ и $DC$, добавляя таким образом в комплекс новую клетку $ABCD$ для каждого пути $ABC$, удовлетворяющим этим условиям. При этом $AB$ и $BC$ будем считать левой и верхней стороной в новой клетке. Вершину $B$ будем называть {\it ядром} нового подклееного подкомплекса.  Полученные вершины $D$ по определению будут уровня $n$.

Будем считать, что на плиточную подстановку наложено дополнительное условие, что в левый верхний угол комплекса любого уровня не входит других ребер, кроме верхней и левой сторон. Из этого условия следует, что из ядра нет ребер внутрь подклееного комплекса, только по его границам.

Пусть $K'_n$ - получившийся комплекс. Он содержит ориентированные плитки (как старые, из $K_n$) так и подклееные. Для каждой ориентированной плитки применим плиточную подстановку, разбивая ее на несколько меньших плиток и добавляя новые вершины. Полученный комплекс будем называть комплексом уровня $n+1$. Добавленные при последнем разбиении вершины будем считать ему принадлежащими и имеющими уровень $n+1$.
\end{definition}

\medskip

{\bf Замечание. } Подклееные плитки при переходе к комплексу большего уровня также подразбиваются. Пути $ABC$ при подклейке вполне могут выбираться и так, что $C$ лежит в подклееной части, а $A$ и $B$ - в плоской.

Кроме того, если вершина $A$ лежит на границе подклееного подкомплекса $K$, и не является его ядром, то она является боковой или краевой, граница $K$ лежит на макроребре, которому принадлежит $A$. Если $A$ - угол $K$, то ровно одно главное выходящее ребро лежит на границе $K$, если $A$ не угол $K$ то оба выходящих главные ребра из $A$ тоже лежат на границе $K$. Неглавные выходящие из $A$ ребра лежат в $K$ только если это ребра в соответствующую подклейку.

\medskip

Докажем некоторую вариацию свойства детерминированности для путей на полученном комплексе. 

Итак, каждый комплекс полученной последовательности состоит из конечного числа плоских подкомплексов (рассмотренных выше). Каждый плоский подкомплекс по определению подклейки двумя своими соседними макроребрами подклеивается к остальной части всего комплекса.

\medskip

Рассмотрим $K_n$ -- плоский комплекс уровня $n$, определенный подстановочным способом, описанным выше. Пусть $ABCD$ — некоторая минимальная плитка на $K_n$. Замену участка пути $ABC$ на участок пути $ADC$ будем называть {\it локальным преобразованием} пути. Если некоторый путь $P’$ получается из пути $P$ с помощью нескольких локальных преобразований, будем считать их {\it эквивалентными}.
 
Путем на комплексе будем называть конечную последовательность 
$V_1$, $E_1$,  $V_2$, $E_2$ \dots $V_{n-1}$, $E_{n-1}$, $V_{n}$, где  $V_i$ - вершины, а $E_j$ - ребра. 
  
Рассмотрим произвольный путь $P$ на $K_n$, концы которого лежат на границе $K_n$. Если для любого такого пути $P$ 
его можно локально преобразовать в путь $P’$, полностью лежащий на границе $K_n$, будем называть $K_n$ {\it плоским регулярным} комплексом.
 
Последовательность плоских регулярных комплексов построена в работе \cite{nilsemigroup1}.
 
\begin{definition} Определение регулярного комплекса.

1. Плоский регулярный комплекс уровня $n$ будем считать регулярным комплексом.

2. Пусть регулярный комплекс $R$ и плоский регулярный комплекс $K_n$ уровня $n$ таковы, что две соседних стороны $K_n$ принадлежат $R$, причем других пересечений у этих комплексов нет. Тогда объединение $R \cup K_n$ также будем считать регулярным комплексом.

\end{definition}

Из определения следует, что регулярный комплекс состоит из объединения конечного числа плоских регулярных подкомплексов. Пусть вершина $A$ лежит на пересечении регулярного подкомплекса $K$ и плоского подкомплекса $K'$, а следующая вершина $B$ лежит внутри $K'$ (не на границе). Тогда ребро из $A$ в $B$ будем называть {\it ребром в подклейку}.

\medskip

\begin{lemma}[Выпрямление пути]  
Пусть путь $P$ принадлежит некоторыму регулярному комплексу $K$, причем концы $P$ принадлежат некоторому плоскому подкомплексу $\hat{K}$. Тогда существует $P'\equiv P$ полностью лежащий в $\hat{K}$.
\end{lemma}

Доказательство. Проведем разметку пути. Если путь проходит по ребру, входящему в некоторую подлкейку, отметим его открывающей скобкой, если наоборот, возвращается из подклейки, поставим закрывающую скобку. Если скобок вообще нет, путь уже полностью лежит в подкомплексе $\hat{K}$ . Проведем доказательство по индукции по общему числу скобок в пути, пусть для $k$ скобок утверждение верное.

\smallskip

Если на пути встречается участок, начинающийся с открывающейся скобки и заканчивающийся закрывающей, то на этом участке путь проходит по некоторому плоскому подкомплексу и можно воспользоваться свойством его регулярности, заменив этот участок на проходящий по границе этого подкомплекса. При этом вместо двух скобок в новом пути будет максимум одна.
Для нового пути можно воспользоваться предположением индукции.

Если указанного выше участка нет, тогда сначала идут все закрывающие скобки, потом все открывающие. Все скобки не могут быть одного типа, так как тогда концы не могли бы быть в одном подкомплексе. Выделим первую открывающую и последнюю закрывающую скобки. Тогда путь представляется в форме $A)XBY(C$, где $X,Y$ -- некоторые вершины, подпуть $A$ содержит только закрывающие скобки, а подпуть $B$ только открывающие. Рассмотрим пути $u$ из $X$ до левого края и $v$ от $Y$ до правого края. Заметим, что $u$ $v$ всегда лишь заходит внутрь подклеек. 

Рассмотрим $K_u$ и $K_v$  -- подкомплексы, в которых лежат, соответственно, $u$ и $v$. $K_v$ и $K_u$ могут пересекаться только по границе. Но концы нашего изначального пути (то есть конечные вершины $u$ и $v$) лежат в одном плоском подкомплексе, что означает, что оба пути $u$ и $v$ начинают и заканчиваются на границе $K_u$ и $K_v$, тогда к ним можно применить предположение индукции и уменьшить число скобок.

\bigskip

В основной части настоящей работы показано как определяются цвета вершин и ребер в плоском комплексе.

По построению, каждая вершина принадлежит ровно одному базовому плоскому подкомплексу, кроме того, она может находится на границе нескольких подклееных подкомплексов. Базовый подкомплекс  вершины $X$ будем обозначать как $K_X$.

\medskip

Пусть имеется путь $P=A_1A_2\dots A_k$ проходящий по некоторому комплексу. Для каждого ребра $A_iA_{i+1}$ определено, является ли это ребро входящим из подклееного комплекса, выходящим в подклееный комплекс или плоским. 
%Будем считать, что два последовательных ребра не являются входящим и выходящим (или наоборот). 
Из определения последовательности комплексов следует, что два последовательных ребра не являются одновременно входящими или выходящими.

%Тогда для каждой вершины $A_i$ не более одного смежного с ней ребра (из пути $P$) уходит в подклейку.

\begin{definition}[Веер и флаг вершины]
 Пусть $B(P)_{A_i}$ -- базовый  комплекс для $A_i$, а $L(P)_{A_i}$ и $R(P)_{A_i}$ -- соответствующие подклееные комплексы (определены, если ребра из $A_i$ уходят в подклейки). Все эти комплексы определяются относительно пути $P$. Упорядоченную тройку $B_{A_i}$, $L_{A_i}$, $R_{A_i}$ (где последние два элемента могут быть пустыми) будем считать {\it ассоциированными } комплексами для $A_i$ и обозначать как $\Omega(P)_{A_i}$.

{\it Веером вершины} $A_i$ будем называть упорядоченную пятерку $\Omega_{A_{i-2}}, \Omega_{A_{i-1}},\Omega_{A_{i}}, \Omega_{A_{i+1}}, \Omega_{A_{i+2}}$, то есть это несколько комплексов являющихся базовыми или подклееными для вершин $A_{i-2},A_{i-1},A_i A_{i+1}, A_{i+2}$, некоторые из которых могут совпадать. 

Пусть вершина $A_i$ зафиксирована. Пусть комплекс $K$ входит в веер $A_i$. Пусть вершина $X$, является начальником $A_i$ в одном из ассоциированных с $A_i$
комплексов. Вместо ее типа, будем рассматривать ее {\it обобщенный тип}: упорядоченный набор типов $X$ относительно всех комплексов в веере $A_i$.

{\it Флагом вершины} $A_i$ будем называть упорядоченный набор цветов $A_i$ относительно всех комплексов, входящих в веер $A_i$.

\end{definition}

{\bf Замечание.} Таким образом, в для пространственного пути, в вершине $A_i$ хранится информация не только о ее цвете в трех ассоциированных комплексах, но и ее полном цвете в ассоциированных комплексах для последовательных вершинах, третья из которых это $A_i$. При этом для начальников $A_i$ используются не обычные плоские типы вершин, а типы относительно всех комплексов в веере $A_i$
\medskip

\begin{definition}[о кодировке пути]

Будем считать, что ребро выхода $AB$ определено относительно того базового плоского подкомплекса, которому принадлежит ребро $AB$.

{Кодировкой пути} $A_1\dots A_k$ будем называть упорядоченную последовательность флагов вершин и ребер входов и выходов вдоль пути.
\end{definition}

\begin{lemma}[ о установлении клетки]
Пусть имеется путь $A_1A_2A_3A_4A_5A_6A_7$, причем если три вершины, идут подряд и лежат в одном плоском подкомплексе, то данный подпуть является регулярным.  Тогда по кодировке пути $A_3A_4A_5$ определяется, являются ли соответствующие вершины тремя вершинами в некоторой клетке комплекса.
\end{lemma}

Доказательство. Если оба ребра из $A_4$ -- плоские, то ситуация сводится к уже рассмотренному плоскому случаю.
Если оба ребра $A_4A_3$ и $A_4A_5$ оба выходят в подклейку, то, очевидно искомой клетки не существует.
Пусть, для определенности, ребро $A_4A_5$ как выходящее из $A_4$, является выходом в подклейку, а $A_4A_3$ - плоское.

Пусть $K_{A_4}$ -- базовый плоский подкомплекс для вершины $A_4$. Допустим, ребро $A_3A_4$ принадлежит $K_{A_4}$. В этом случае $A_4$ не может быть ядром 
$K_{A_4}$.  Если это ребро, как входящее в $A_4$, является неглавным, то искомой клетки не существует, так как по границе $K_{A_4}$ проходит макроребро на котором лежит $A_4$.

Пусть ребро $A_3A_4$ как входящее в $A_4$ является главным. Тогда $A_3A_4$ лежит на границе $K_{A_4}$, как и весь путь $A_3A_4A_5$, тогда данный подпуть является регулярным и мы получаем искомое.

\medskip

\begin{lemma}[о детерминированности вдоль пути]
Пусть имеется путь $A_1A_2A_3A_4A_5A_6A_7$ и выполнены условия:

1) $A_3$, $A_4$, $A_5$ являются вершинами в некоторой клетке $T$;

2) Путь $A_3A_4A_5$ является регулярным относительно плоского подкомплекса, содержащего клетку $T$.

Пусть $X$ -- четвертая вершина этой клетки. Тогда  по известной кодировке пути $A_2A_3A_4A_5A_6$ однозначно определяется кодировка пути $A_2A_3XA_5A_6$. 
\end{lemma}

Доказательство. 
случай 1. Оба смежных с $A_4$ ребра являются плоскими. В этом случае кусок пути $A_3A_4A_5$ является плоским. К нему применима теорема о слабой детерминированности в плоском случае, таким образом по цветам $A_3$, $A_4$, $A_5$ в базовом комплексе однозначно определяется цвет $X$ в базовом комплексе.
Определим остальные компоненты флага $X$. 
 Допустим, из $A_5$ следующее ребро пути уходит в подклееный комплекс $K_{A_5A_6}$.  Фактически, требуется определить цвет $X$ относительно $K_{A_5A_6}$. Заметим, что $X$ и $A_5$ лежат на макроребре, общем для комплексов $K_{A_5A_6}$ и базовом, где лежит путь $A_3A_4A_5$. Тогда клетка $T=XA_3A_4A_5$ примыкает к этому макроребру вершинами $X$ и $A_5$. Рассмотрим минимальный плоский подкомплекс, содержащий клетку $T$. Тогда либо $A_3$ либо $A_4$ лежат внутри него, и $X$ является начальником соответствующей вершины, следующим или предыдущим по часовой стрелке в зависимости от ориентации пути $A_3A_4A_5$ внутри $T$. Таким образом, тип $X$ относительно $K_{A_5A_6}$ определяется. Начальники $X$ относительно $K_{A_5A_6}$ совпадают с начальниками $A_5$.
  
 Таким образом, полностью определен флаг вершины $X$.

\begin{figure}[hbtp]
\centering
\includegraphics[width=0.3\textwidth]{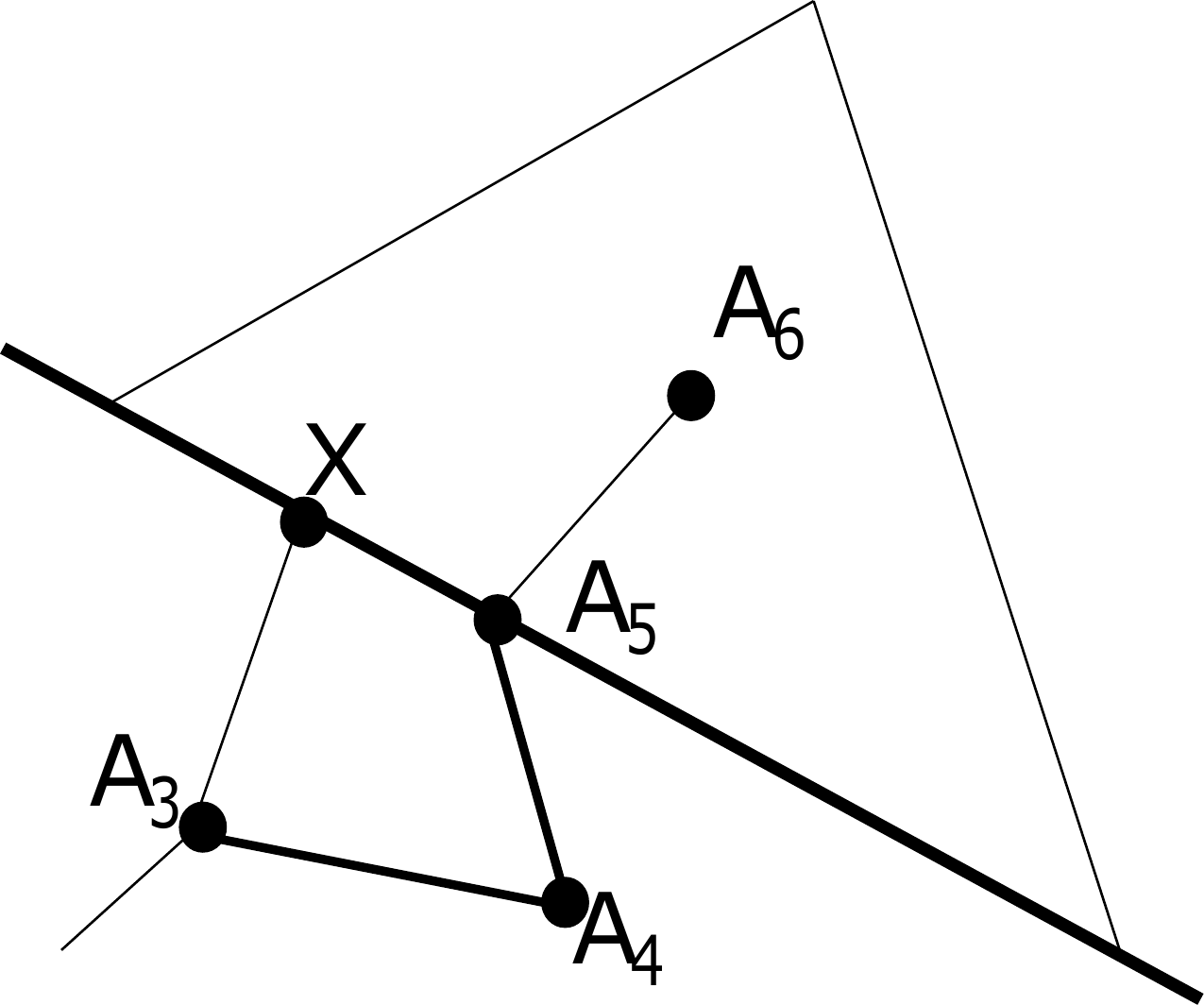}
\caption{Восстановление $X$, лежащей на общем макроребре с подклееным комплексом}
\label{fig:flags}
\end{figure}

\medskip

случай 2. Одно из ребер, смежных с $A_4$ является выходящим в подклейку. Будем считать, что это ребро $A_4A_5$. Рассмотрим комплекс $K$, являющийся базовым по отношению к вершине $A_5$. Тогда $A_3$ и $A_4$ лежат на краю $K$ (иначе клетки $T$ не существует). Рассмотрим путь $A_3A_4A_5$ как плоский в $K$. Тогда по теореме о плоской детерминированности однозначно определяется цвет $X$ относительно $K$. Заметим, что, веер вершин $X$ и $A_4$ совпадает(помимо самих вершин $X$ и $A_4$).  Таким образом, флаг вершины $X$ полностью определяется, аналогично плоскому случаю.

\medskip

Лемма выше является пространственным аналогом теоремы о слабой детерминированности. Она используется для введения полугрупповых определяющих соотношений для
уже построенной последовательности подстановочных комплексов. После введения данных соотношений появляется возможность преобразовывать слова, отвечающие кодировкам путей. Преобразования слов отвечают переходам к эквивалентным путям.

\end{document}